\documentclass[11pt]{amsart}
\usepackage{amsmath,amssymb,amsfonts,times,hyperref}

\title
{Optimal Distortion Embeddings of Distance Regular Graphs into
Euclidean Spaces}

\author{Frank Vallentin}
\address{F. Vallentin, Centrum voor Wiskunde en Informatica (CWI), Kruislaan 413,
1098 SJ Amsterdam, The Netherlands} 
\email{f.vallentin@cwi.nl}

\thanks{ The author was supported by the Netherlands Organization for
Scientific Research under grant NWO 639.032.203 and by the Deutsche
Forschungsgemeinschaft (DFG) under grant SCHU 1503/4-1.  Part of this
work was done when the author was supported by the Edmund Landau
Center for Research in Mathematical Analysis and Related Areas,
Einstein Institute of Mathematics, The Hebrew University Jerusalem,
Israel, sponsored by the Minerva Foundation (Germany).}

\date{24th May 2007}


\newcommand{\R}{\mathbb R}
\newcommand{\Z}{\mathbb Z}

\newcommand{\MA}{\mathcal A}

\DeclareMathOperator{\Aut}{Aut}
\DeclareMathOperator{\diam}{diam}
\DeclareMathOperator{\Sym}{Sym}
\DeclareMathOperator{\trace}{trace}
\DeclareMathOperator{\card}{card}

\newtheorem{definition}{Definition}[section]
\newtheorem{lemma}[definition]{Lemma}
\newtheorem{theorem}[definition]{Theorem}

\newtheorem{remark}[definition]{Remark}
\newtheorem{example}[definition]{Example}

\begin{document}

\begin{abstract}
In this paper we give a lower bound for the least distortion embedding
of a distance regular graph into Euclidean space. We use the lower
bound for finding the least distortion for Hamming graphs,
Johnson graphs, and all strongly regular graphs. Our technique
involves semidefinite programming and exploiting the algebra structure
of the optimization problem so that the question of finding a lower
bound of the least distortion is reduced to an analytic question about
orthogonal polynomials.
\end{abstract}

\maketitle

\section{Introduction}
\label{sec:introduction}

By $\R^n$ we denote the Euclidean space of column vectors $x = (x_1,
\ldots, x_n)^t$ with standard inner product $x \cdot y= x_1 y_1 +
\cdots + x_n y_n$ and corresponding norm $\|x\| = \sqrt{x \cdot
x}$. Let $(X,d)$ be a finite metric space with $n$ elements. We say
that an embedding $\varrho : X \to \R^n$ into Euclidean space has
\emph{distortion} $D$ if for all $x,y \in X$ the inequalities
\[
d(x,y) \leq \|\varrho(x) - \varrho(y)\| \leq D d(x,y)
\]
hold. 

By $c_2(X,d)$ we denote the \emph{least distortion} for which
$(X,d)$ can be embedded into $\R^n$ and say that an embedding of $(X,
d)$ is \emph{optimal} if it has distortion $c_2(X,d)$.

In \cite{bourgain-1985} Bourgain showed that $c_2(X,d) = O(\log n)$
and in \cite{llr-1995} Linial, London and Rabinovich proved that this
bound is tight. In the last years embeddability questions, especially
of finite graphs where the metric is given by the shortest path
metric, were studied by theoretical computer scientists. For example
they were used to design approximation algorithms (see e.g.\
\cite{linial-2002}, \cite{indyk-2001} and \cite{matousek-2002},
Chapter 15).

Despite this interest for only very few graphs the exact least
distortion and a least distortion embedding is explicitly known. The
list only includes unit cubes (due to Enflo, see \cite{enflo-1969}),
cycles, and strong graph product of cycles (due to Linial and Magen,
see \cite{lm-2000}). Extending work of Linial and Magen we give a
lower bound for the least distortion of distance regular graphs. It
turns out that the bound is tight in many examples and \emph{we
conjecture that it is always tight}. We compute least distortions for
the following important examples: Hamming graphs (which include the
cube), Johnson graphs, and all strongly regular graphs.

This paper is organized as follows. In Section~\ref{sec:results} we
give the necessary definitions and state our results. In
Section~\ref{sec:proof} we prove the lower bound and in
Section~\ref{sec:examples} we work out the three cases.

\section{Statement of Results}
\label{sec:results}

Before we formulate our results we recall some definitions and results
of the theory of distance regular graphs. For a comprehensive
treatment we refer to \cite{bi-1984} and \cite{bcn-1989}.

Let $G = (V, E)$ be an \emph{undirected graph} given by a finite set
$V$ of \emph{vertices} and a subset $E \subseteq \binom{V}{2}$ of
two-element subsets of $V$ called \emph{edges}. By $d : V \times V
\to \Z_{\geq 0} \cup \{\infty\}$ we denote the length of a shortest
path connecting two vertices $x$ and $y$ in $G$ where we set $d(x, y)
= \infty$ whenever there is no connection at all. The
\emph{diameter} of $G$ is $\diam G = \max_{x, y \in V} d(x, y)$. A
\emph{connected} graph $G$, that is a graph with finite diameter,
gives a finite metric space $(V, d)$. In this situation we write for
the least distortion $c_2(G)$ instead of $c_2(V, d)$.

A connected graph $G$ is called \emph{distance regular} if there are
constants $a_i, b_i, c_i$ where $i \in \{0, \ldots, \diam G\}$ so that
the following holds: For every pair of vertices $x, y \in V$ with
$d(x, y) = i$ we have
\begin{equation}
\begin{array}{lll}
a_i & = & \card(\{ z \in V : \mbox{$d(x, z) = 1$ and $d(z, y) = i$}\}),\\
b_i & = & \card(\{ z \in V : \mbox{$d(x, z) = 1$ and $d(z, y) = i + 1$}\}),\\
c_i & = & \card(\{ z \in V : \mbox{$d(x, z) = 1$ and $d(z, y) = i - 1$}\}).
\end{array}
\end{equation}
The number
\begin{equation}
k_i = \card(\{ y \in V : d(x, y) = i\}), 
\end{equation}
is called the \emph{$i$-th degree} of $G$. It is independent
of $x$.

The following three families are important examples of distance
regular graphs. We will find their least distortions in
Section~\ref{sec:examples}.

\begin{example} [Hamming Graphs]
\label{ex:hamming}
Let $X$ be a finite set of cardinality $q \geq 2$. The vertex set of
the \emph{Hamming graph} $H(q, n)$ is $X^n$, the set of all vectors
of length $n$. Two vertices $x, y \in X^n$ are adjacent if $x$ and $y$
differ in exactly one coordinate. The shortest path metric of $H(q,
n)$ coincides with the Hamming distance. The diameter of $H(q, n)$ is
$n$.
\end{example}

\begin{example}[Johnson Graphs] 
Let $V$ be a set of size $v$ and $n$ be an integer with $v \geq
2n$. The vertex set of the \emph{Johnson graph} $J(v, n)$ is the set
$\binom{V}{n}$ of all $n$-element subsets of $V$. Two vertices $x, y$
of $J(v, n)$ are adjacent if the intersection $x \cap y$ has
cardinality $n - 1$. The diameter of $J(v, n)$ is $n$.
\end{example}

\begin{example}[Strongly Regular Graphs]
A \emph{strongly regular graph} with parameters $(\nu, k, \lambda,
\mu)$ is a graph with $\nu$ vertices where every vertex is adjacent to
$k$ vertices, where every pair of adjacent vertices has precisely
$\lambda$ common neighbors, and where every pair of nonadjacent
vertices has precisely $\mu$ common neighbors. If a strongly regular
graph has diameter $2$, then it is a distance regular graph with $k_1
= k$, $a_1 = \lambda$, $c_2 = \mu$. Otherwise it is a disjoint union
of equal-sized complete graphs.
\end{example}

For $i \in \{0, \dots, \diam G\}$ we define the \emph{$i$-th
adjacency matrix} $A_i \in \{0, 1\}^{V \times V}$ component wise by
$(A_i)_{xy} = 1$ whenever $d(x, y) = i$ and $(A_i)_{xy} = 0$
otherwise. We have the following relation between the adjacency
matrices
\begin{equation}
\label{eq:adjacencyrelation}
A_1 A_i = c_{i + 1} A_{i + 1} + a_i A_i + b_{i - 1} A_{i - 1}.
\end{equation}
Hence we can write $A_i = v_i(A_1)$ for univariate polynomials $v_i$
of degree $i$. By $\theta_0 > \ldots > \theta_{\diam G}$ we denote the
different eigenvalues of the matrix $A_1$. Notice that $v_i(\theta_0)
= k_i$ and that $v_i(\theta_0)$ is the largest eigenvalue of $A_i$.

\medskip

\noindent Now we can state our principal theorem.

\begin{theorem}
\label{th:main}
Let $G$ be a distance regular graph with $n = \diam G$. Then,
\begin{equation}
c_2(G)^2 \geq \frac{n^2 v_n(\theta_0)}{v_1(\theta_0)} \min_{j \in \{1,
\ldots, n\}} \Big\{\frac{v_1(\theta_0) - v_1(\theta_j)}{v_n(\theta_0)
- v_n(\theta_j)} \Big\}.
\end{equation}
\end{theorem}

\noindent We prove this theorem in Section~\ref{sec:proof}. The proof
is based on the following observations. In general one can compute a
least distortion embedding by solving a semidefinite programming
problem. Using the commutativity of the algebra spanned by the
adjacency matrices one can simplify the semidefinite programming
program considerably (even to a linear program, see
e.g. \cite{gr-1999}). Then, using duality theory of semidefinite
programming one gets a lower bound for the least distortion.

\medskip

In Section~\ref{sec:examples} we apply this theorem to the distance
regular graphs we introduced above to get their least distortions. The
following theorem summarizes the results.

\begin{theorem}
\label{th:examples}
\noindent 
\begin{enumerate}
\item[(a)] For the Hamming graph $H(q, n)$ we have 
\[
c_2(H(q, n)) = \sqrt{n}.
\]
\item[(b)] For the Johnson graph $J(v, n)$ we have
\[
c_2(J(v, n)) = \sqrt{n}.
\]
\item[(c)] For a strongly regular graph $G$ of diameter $2$ with parameters
$(\nu, k, \lambda, \mu)$ we have
\[
c_2(G) = \sqrt{\frac{4(\nu-k-1)(k-r)}{k(\nu-k+r)}}, 
\]
where $r = \frac{1}{2}\big(\lambda - \mu + \sqrt{\nu}\big)$.
\end{enumerate}
\end{theorem}

\section{Proof of Theorem~\ref{th:main}}
\label{sec:proof}

Linial, London and Rabinovich \cite{llr-1995} were the first who
noticed that finding a least distortion embedding of a finite metric
space $(X, d)$ into Euclidean space can be expressed as a semidefinite
programming problem:
\begin{equation}
\label{eq:sdp}
\begin{array}{ll}
  \mbox{\textbf{minimize}} & C\\[1ex]
\mbox{\textbf{subject to}} & \mbox{$Q = (q_{xy}) \in \R^{X \times X}$ is
positive semidefinite},\\[1ex]
& \mbox{$d(x, y)^2 \leq q_{xx} - 2q_{xy} + q_{yy} \leq C d(x, y)^2$ for all $x,y \in X$.}
\end{array}
\end{equation}
Here $Q$ is the \emph{Gram matrix} of an embedding $\varrho : X \to
\R^{n}$ defined entry wise by $q_{xy} = \varrho(x) \cdot
\varrho(y)$. Note that $Q$ defines the embedding $\varrho$ uniquely up
to orthogonal transformations. The minimum $C$ of the semidefinite
programming problem~(\ref{eq:sdp}) equals $c_2(X, d)^2$.

Semidefinite programming problems are convex minimization problems and
they can be solved efficiently in polynomial time in the sense that
one can approximate an optimal solution to any fixed precision (see
the survey \cite{vb-1996}). Furthermore, semidefinite programming
problems respect the symmetries of the instances. Hence, there is a
least distortion embedding of a distance regular graph which inherits
the symmetries of the graph. Now we make this statement precise. For
this we start with a definition.

\begin{definition}
\label{def:faithful}
Let $(X, d)$ be a finite metric space. We say that an embedding
$\varrho : X \to \R^{n}$ into Euclidean space is \emph{faithful} if for
every two pairs $(x, y)$ and $(x', y') \in X \times X$ we have
\begin{equation}
d(x, y) = d(x', y') \Longrightarrow \|\varrho(x) - \varrho(y)\| =
\|\varrho(x') - \varrho(y')\|.
\end{equation}
\end{definition}

\begin{lemma}
\label{lem:faithful}
Let $G = (V, E)$ be a distance regular graph. Then, there exists a
faithful embedding of $G$ into Euclidean space with minimal
distortion.
\end{lemma}

\begin{proof}
Let $Q \in \R^{V \times V}$ be the Gram matrix of an embedding
$\varrho : V \to \R^{n}$. We denote the entries of $Q$ by $q_{xy} =
\varrho(x) \cdot \varrho(y)$. Suppose that $\varrho$ has distortion
$D$ so that we have the inequality
\begin{equation}
\label{eq:distortion}
d(x, y)^2 \leq q_{xx} - 2q_{xy} + q_{yy} \leq D^2 d(x, y)^2
\end{equation}
for all $x, y \in V$. 

Because of \eqref{eq:adjacencyrelation} the algebra $\MA$ generated by
the adjacency matrices $A_i$ is commutative. The algebra $\MA$ is
called the \emph{Bose-Mesner algebra} of $G$ and it has basis $A_i$
with $i = 0, \ldots, \diam G$. It is equipped with the inner product
$\langle A, B \rangle = \trace(A^t B)$.

Now we show that the orthogonal projection $\bar{Q}$ of $Q$ onto $\MA$
is a Gram matrix of a faithful embedding having distortion $D$. 

First we argue that $\bar{Q}$ is positive semidefinite.
Because $\MA$ is commutative the adjacency matrices $A_i$ have a
common basis of eigenvectors. Decompose the space $\R^V$
into an orthogonal direct sum of maximal common eigenspaces:
\begin{equation}
\R^V = V_0 \perp V_1 \perp \ldots \perp V_{\diam G}.
\end{equation}
Then, the matrices of the orthogonal projection $E_i : \R^V \to V_i$
form a basis of $\MA$. Since they are positive semidefinite we have
$\langle Q, E_i \rangle \geq 0$. Hence the orthogonal projection
\begin{equation}
\bar{Q} = \sum_{i = 0}^{\diam G} \frac{\langle Q,
E_i\rangle}{\langle E_i, E_i \rangle} E_i
\end{equation}
is positive semidefinite.

To show that $\bar{Q}$ is faithful and satisfies the desired
inequalities we use the representation
\begin{equation}
\bar{Q} = \sum_{i = 0}^{\diam G} \frac{\langle Q, A_i\rangle}{\langle
A_i, A_i \rangle} A_i.
\end{equation}
Notice here that the adjacency matrices form an orthogonal basis of
$\MA$. Let $x, y \in V$ be two vertices at distance $d = d(x,y)$. For
the entry $\bar{q}_{xy}$ of $\bar{Q}$ we have
\begin{equation}
\label{eq:faithful}
\bar{q}_{xy} = \sum_{i=0}^{\diam G} \frac{\langle Q,
A_i\rangle}{\langle A_i, A_i \rangle} (A_i)_{xy} = \frac{\langle Q,
A_d\rangle}{\langle A_d, A_d\rangle} = \frac{1}{\card(M_d)}\sum_{(x',
y') \in M_d} q_{x'y'},
\end{equation}
where $M_d = \{(x,y) \in V \times V: d(x,y) = d\}$. From
\eqref{eq:faithful} it follows immediately that the embedding
$\bar{\varrho}$ given by $\bar{Q}$ is faithful. Furthermore we
obviously have
\begin{equation}
d(x, y)^2 = \displaystyle \frac{1}{\card(M_d)}
\sum_{(x',y') \in M_d} d(x', y')^2.
\end{equation}
Applying this to (\ref{eq:distortion}) and using the definition of
$\bar{Q}$ gives
\begin{equation}
\begin{array}{rcl}
d(x, y)^2 
& \leq & \displaystyle \frac{1}{\card(M_d)} \sum_{(x',y') \in M_d} \big(q_{x'x'} - 2q_{x'y'} + q_{y'y'}\big)\\[3ex]
& = & \displaystyle \bar{q}_{xx} - 2\bar{q}_{xy} + \bar{q}_{yy} \\[1ex]
& \leq & \displaystyle \frac{D^2}{\card(M_d)} \sum_{(x',y') \in M_d} d(x', y')^2\\[3ex]
& = & \displaystyle D^2 d(x, y)^2,
\end{array}
\end{equation}
hence the embedding given by $\bar{Q}$ has distortion $D$.
\end{proof}

\begin{remark}
If the graph $G$ is distance transitive, then one can partially
simplify the proof of Lemma~\ref{lem:faithful}: The \emph{automorphism
group} $\Aut(G)$ is the set of permutations $\sigma \in \Sym(V)$ with
$\text{$\{x, y\} \in E$ if and only if $\{\sigma(x), \sigma(y)\} \in
E$}$, and we say that $G$ is \emph{distance transitive} if for every
pair of vertex pairs $(x, y), (x', y')$ with $d(x, y) = d(x', y')$
there exists $\sigma \in \Aut(G)$ so that $(\sigma(x), \sigma(y)) = (x',
y')$. Then, the orthogonal projection $\bar{Q}$ is simply the
symmetrization of $Q$, i.e.\
\begin{equation}
\bar{Q} =
\frac{1}{|\Aut(G)|}\sum_{\sigma \in \Aut(G)} (q_{\sigma(x),
\sigma(y)}),
\end{equation}
and $\bar{Q}$ is positive semidefinite because it is the sum of
positive semidefinite matrices.
\end{remark}

Using duality theory of semidefinite programming Linial, London and
Rabinovich \cite{llr-1995} and Linial and Magen \cite{lm-2000} gave the
following characterization of the least possible distortion for a
finite metric space.

\begin{lemma}
\label{lem:criterion}
Let $(X, d)$ be a finite metric space. 

\begin{enumerate}
\item[(a)]
The least distortion of an embedding of $(X, d)$ into Euclidean space
is given by
\begin{equation}
\label{eq:criterion}
c_2(X, d)^2 = \max_Q \frac{\sum\limits_{\{(x, y) : q_{xy} > 0\}} d(x,
y)^2 q_{xy}}{\sum\limits_{\{(x, y) : q_{xy} < 0\}} d(x, y)^2 (-q_{xy})},
\end{equation}
where the maximum is taken among all positive semidefinite matrices
$Q$ in which all row sums vanish. (Note that the quotient is invariant
under scaling of $Q$ with positive reals.)

\item[(b)] Let $\varrho$ be an embedding of $(X, d)$ into Euclidean
space having minimal distortion $c_2(X, d)$. For a matrix attaining
the maximum in (\ref{eq:criterion}) and for a pair $(x, y) \in X
\times X$ we have $q_{xy} > 0$ only for the \emph{most contracted
pairs} $(x, y)$, that is for $(x, y)$ the fraction $\|\varrho(x) -
\varrho(y)\| / d(x, y)$ is minimal among all pairs in $X \times X$, we
have $q_{xy} < 0$ only for the \emph{most expanded pairs} $(x, y)$,
that is for $(x, y)$ the fraction $\|\varrho(x) - \varrho(y)\| / d(x,
y)$ is maximal among all pairs in $X \times X$, and $q_{xy} = 0$ for
all other pairs.
\end{enumerate}
\end{lemma}

\begin{proof}
See \cite[Corollary 3.5]{llr-1995} and \cite[Claim 1.4]{lm-2000}.
\end{proof}

\begin{remark}
For the embedding of finite metric spaces given by the shortest path
metric of a graph, Linial and Magen showed (\cite[Claim 2.2]{lm-2000})
that most expanded pairs are always adjacent vertices.
\end{remark}

Now we finish the proof of Theorem~\ref{th:main}. Let $G$ be a
distance regular graph and let $\varrho$ be an embedding of $G$ into
Euclidean space with minimal distortion~$c_2(G)$. By
Lemma~\ref{lem:faithful} we can assume that $\varrho$ is
faithful. Hence, by the previous remark, all pairs $(x, y)$ with $d(x,
y) = 1$ are most expanded, and there is an index $i \in \{2, \ldots,
\diam G\}$ so that all pairs $(x, y)$ with $d(x, y) = i$ are most
contracted.

For proving a lower bound on the distortion of $\varrho$ we suppose
that $i = n$, where $n = \diam G$. So the lower bound can only be
tight when the most contracted pairs are at distance $n$.

We define
\begin{equation}
Q_{\alpha} = (k_1 - \alpha k_{n}) A_0 - A_1 + \alpha A_{n},\;\;\; \alpha \in \R.
\end{equation}
When $Q_{\alpha}$ is positive semidefinite, then $Q_{\alpha}$ satisfies
the assumption of Lemma~\ref{lem:criterion} (a). Hence,
\begin{equation}
c_2(G)^2  \geq  \Big\{\frac{k_{n} n^2
\alpha }{k_1} :
\mbox{$Q_{\alpha}$ is positive semidefinite}
\Big\}.
\end{equation}
In order to maximize $\frac{k_{n} n^2 \alpha }{k_1}$ we have to
maximize $\alpha$ so that $Q_{\alpha}$ is positive semidefinite.
Recall that the adjacency matrices have a common system of
eigenvectors.  Let $x_j$ be a common eigenvector of the adjacency
matrices which is an eigenvector of the eigenvalue $\theta_j$ of
$A_1$. Then, $A_1 x_j = \theta_j x_j$, and
\begin{equation}
 Q_{\alpha}x_j
= (k_1 - \alpha k_{n} - \theta_j +  \alpha v_{n}(\theta_j)) x_j,
\end{equation}
and the matrix $Q_{\alpha}$ is positive semidefinite if and only if
\begin{equation}
k_1 - \alpha k_{n} - \theta_j + \alpha v_{n}(\theta_j) \geq 0, \quad \mbox{for
all $j \in \{0, \ldots, n\}$.}
\end{equation}
The largest eigenvalue of the adjacency matrix of a $k$-regular graph
is exactly $k$. So, $k_{n} - v_{n}(\theta_j)$ is positive for $j \in
\{1, \ldots, n\}$ and $k_{n} - v_{n}(\theta_j) = 0$ for $j =
0$. Hence,
\begin{equation}
\begin{array}{rcl}
\alpha & = & \displaystyle\min_{j \in \{1,
\ldots, n\}} \frac{k_1 - \theta_j}{k_n - v_n(\theta_j)}\\[2ex]
& = & \displaystyle\min_{j \in \{1,
\ldots, n\}} \frac{v_1(\theta_0) - v_1(\theta_j)}{v_n(\theta_0) - v_n(\theta_j)},\\
\end{array}
\end{equation}
which yields the statement of the theorem.

\section{Examples}
\label{sec:examples}

\subsection{Hamming Graphs}
\label{ssec:hamming}

Now we show using Theorem~\ref{th:main} that the optimal distortion of
the Hamming graph $H(q, n)$ is $\sqrt{n}$ and we give an embedding of
$H(q, n)$ into Euclidean space having this distortion.

We use the notation we introduced in Section~\ref{sec:results}. The
eigenvalues of the $i$-th adjacency matrix of $H(q, n)$ are well-known
(see for example \cite[Chapter 3.2]{bi-1984}). They are $v_i(\theta_j) =
K_i(j)$ where $j \in \{0, \ldots, n\}$ and where $K_i$ is the
\emph{$i$-th Krawtchouk polynomial}
\begin{equation}
K_i(u) = \sum_{t = 0}^i (-q)^t (q - 1)^{(i - t)} \binom{n - t}{i - t}.
\binom{u}{t}.
\end{equation}
In particular we have
\begin{equation}
k_i = K_i(0) = \binom{n}{i} (q - 1)^i,
\end{equation}
\begin{equation}
\theta_j = K_1(j) = n(q - 1) - qj,
\end{equation}
\begin{equation}
v_n(\theta_j) = (-1)^j (q - 1)^{n - j}.
\end{equation}
Let us determine the value of $\alpha = \min_{j \in \{1, \ldots,
n\}} \frac{k_1 - \theta_j}{k_n - v_n(\theta_j)}$. The minimum is
attained for $j = 1$ so that we have
\begin{equation}
\alpha = \frac{k_1 -\theta_1}{k_n - v_n(\theta_1)} = \frac{1}{(q - 1)^{n-1}},
\end{equation}
since for $j = 2, \ldots, n$ the inequality
\begin{equation}
\frac{k_1 -\theta_j}{k_n - v_n(\theta_j)} = \frac{qj}{(q - 1)^n -
(-1)^n (q - 1)^{n - j}} \geq \frac{1}{(q-1)^{n-1}}
\end{equation}
holds true. Hence by Theorem~\ref{th:main} the distortion of an
optimal embedding is bounded by
\begin{equation}
c_2(H(q, n))^2 \geq \frac{n^2 \alpha k_n}{k_1} = n.
\end{equation}
We have equality since the embedding $\varrho$ we define below
has distortion $\sqrt{n}$. Let $X^n$ be the vertex set of $H(q,
n)$. With $e_x \in \R^X$ denote the standard unit vector defined
component wise by $(e_x)_y = 1$ if $x = y$ and $(e_x)_y = 0$
otherwise. For a vertex $(x_1, \ldots, x_n) \in X^n$ in $H(q, n)$ set
\begin{equation}
\varrho(x_1, \ldots, x_n) = \sqrt{n/2}(e_{x_1}, \ldots, e_{x_n})^t
\in (\R^X)^n.
\end{equation}
If $d((x_1, \ldots, x_n), (y_1, \ldots, y_n)) = i$, then
$\|\varrho(x_1, \ldots, x_n) - \varrho(y_1, \ldots, y_n)\| =
\sqrt{ni}$ and we have the desired inequalities
\begin{equation}
d(x,y)^2 = i^2 \leq \|\varrho(x) - \varrho(y)\|^2 = ni \leq n d(x,y)^2
= n i^2,
\end{equation}
where we abbreviate $(x_1, \ldots, x_n)$ and $(y_1, \ldots, y_n)$ by
$x$ and $y$. The image of this embedding forms the vertex set of the
direct product, taken $n$ times, of a regular simplex with $q$
vertices.

\begin{remark}
In particular this implies the classical result of Enflo
\cite{enflo-1969} that the least distortion embedding of the
$n$-dimensional unit cube $H(2, n)$ is $\sqrt{n}$. Enflo's proof uses
inductive and combinatorial arguments and does not easily generalize
to different finite metric spaces. Linial and Magen \cite[Theorem 2.4]{lm-2000}
give another proof of Enflo's theorem which is in a sense an ad-hoc
variant of our proof.
\end{remark}

\subsection{Johnson Graphs}
\label{ssec:johnson}

Here we show that the optimal distortion of the Johnson graph $J(v,
n)$ is $\sqrt{n}$ and we give an embedding of $J(v, n)$ into Euclidean
space having this distortion.

The eigenvalues of the $i$-th adjacency matrix of $J(v, n)$ are
well-known (see for example \cite[Chapter 3.2]{bi-1984}). They are
$v_i(\theta_j) = E_i(j)$ where $E_i$ is the \emph{$i$-th Eberlein
polynomial} (or \emph{dual Hahn polynomial})
\begin{equation}
E_i(u) = \sum_{t=0}^i (-1)^t \binom{u}{t} \binom{n - u}{i - t}
\binom{v - n- u}{i - t}.
\end{equation}
In particular we have 
\begin{equation}
k_i = E_i(0) = \binom{n}{i} \binom{v - n}{i},
\end{equation}
\begin{equation}
\theta_j = E_1(j) = j^2 - (v + 1) j + n(v - n).
\end{equation}

Let us determine the value of $\alpha$. We have
\begin{equation}
\label{eq:alphanj}
\alpha = \min_{j = 1, \ldots, n}
\frac{k_1 - \theta_j}{k_n - v_n(\theta_j)} =
\min_{j = 1, \ldots, n}
\frac{(v + 1)j - j^2}{\binom{v - n}{n} - (-1)^j \binom{v - n - j}{n - j}}.
\end{equation}
We shall show that the minimum is attained for $j = 1$ so that
\begin{equation}
\label{eq:alphan}
\alpha =  \frac{v}{\binom{v - n}{n} + \binom{v - n - 1}{n - 1}}.
\end{equation}
We compare the numerator of the right hand side of (\ref{eq:alphan})
with the one of (\ref{eq:alphanj}). This gives the following inequality
which holds true for all $j$ in the interval $[1, m]$
\begin{equation}
v \leq (v + 1) j - j^2.
\end{equation}
We compare the denominators getting the inequality
\begin{equation}
\binom{v - n}{n} - (-1)^j \binom{v - n - j}{n - j} \leq \binom{v - n}{n} +
\binom{v - n - 1}{n - 1},
\end{equation}
which holds because $\binom{v - n - 1}{n - 1} = \binom{v - n - j}{n -
j} \prod_{t = 1}^{j - 1} \frac{v - n - t}{n - t}$ and $v - n - t \geq
n - t$ since $v \geq 2n$. Altogether this shows that the value
$\alpha$ is the one stated in (\ref{eq:alphan}). Hence the squared
distortion of an embedding is at least $n$.  We have equality since
the embedding $\varrho$ described below has distortion $\sqrt{n}$.

Let $\binom{V}{n}$ be the vertex set of $J(v, n)$. With $e_v \in \R^V$
denote the standard unit vector as in the last section.  For a
$n$-element subset $X \subseteq V$ define the embedding $\varrho(X) =
\sqrt{n} \sum_{x \in X} e_x$. If two $n$-element subsets $X$, $Y$ have
distance $i$ in $J(v, n)$, then $\|\varrho(X) - \varrho(Y)\| =
\sqrt{ni}$. Hence, the distortion of $\varrho$ is $\sqrt{n}$. The
image of this embedding forms the vertex set of the $n$-hypersimplex
in dimension $v$.

\subsection{Strongly Regular Graphs}
\label{ssec:stronglyregular}

In this section we will show that the optimal distortion of a strongly
regular graph $G = (V, E)$ of diameter $2$ with parameters $(\nu, k,
\lambda, \mu)$ is
$\big(\frac{4(\nu-k-1)(k-r)}{k(\nu-k+r)}\big)^{1/2}$, where $r =
\frac{1}{2}\big(\lambda - \mu + \sqrt{\nu}\big)$. In the following we
shall make use of \cite[Theorem 1.3.1]{bcn-1989} where fundamental
facts about the parameters $\nu, k, \lambda, \mu$ are provided.

The eigenvalues of the first adjacency matrix $A_1$ are
\begin{equation}
k,\;\; r = \frac{1}{2}\big(\lambda - \mu + \sqrt{\nu}\big),\;\; s = \frac{1}{2}\big(\lambda -
\mu - \sqrt{\nu}\big).
\end{equation}
We have
\begin{equation}
A_1^2 = kA_0 + \lambda A_1 + \mu A_2,
\end{equation}
and hence
\begin{equation}
v_2(u) = \frac{1}{\mu}\big(u^2 - \lambda u - k\big).
\end{equation}
Using the identities $\lambda = \mu + r + s$ and $rs = \mu - k$ we
compute $v_2(r) = -r - 1$ and $v_2(s) = -s - 1$.  Because $r \geq 0$
and $s \leq -1$ we have the inequality
\begin{equation}
\frac{k-r}{(\nu-k-1) - (-1-r)} \leq \frac{k-s}{(\nu-k-1) - (-1-s)}
\end{equation}
Now Theorem~\ref{th:main} gives the lower bound
\begin{equation}
c_2(G)^2 \geq \frac{4(\nu-k-1)(k-r)}{k(\nu-k+r)}.
\end{equation}

By reviewing the proof of Theorem~\ref{th:main} for the case of
distance regular graphs with diameter $2$, i.e.\ for connected
strongly regular graphs, one sees that Theorem~\ref{th:main} is tight
in these cases. The reason for this is that in a faithful embedding
all the most contracted pairs are pairs of vertices which are not
adjacent. So this case is especially convenient since we do not have
to construct an embedding to upper bound the least distortion.

\section*{Acknowledgments}

The author thanks Nati Linial for proposing the project of computing
the least distortion embeddings of highly symmetric graphs. He thanks
Dion Gijswijt, Gil Kalai, Achill Sch\"urmann, Uli Wagner, and the
anonymous referees for useful comments and suggestions.

\bibliographystyle{[15]}

\begin{thebibliography}{BCN89}

\bibitem{bi-1984}
E. Bannai and T. Ito, \emph{Algebraic combinatorics I: Association schemes},
  Benjamin/Cummings, Menlo Park, California, 1984.

\bibitem{bcn-1989} 
A.E. Brouwer, A.M. Cohen and A. Neumaier, \emph{Distance-regular
graphs}, Springer-Verlag, Berlin, 1989.

\bibitem{bourgain-1985}
J. Bourgain, \emph{On Lipschitz embeddings of finite metric spaces in
 Hilbert spaces}, Israel J. Math. \textbf{52} (1985), 46--52.

\bibitem{enflo-1969} 
P. Enflo, \emph{On the nonexistence of uniform homeomorphisms between
{$L\sb{p}$}-spaces}, Ark. Mat. \textbf{8} (1969), 103--105.

\bibitem{gr-1999} 
M.X. Goemans and F. Rendl, \emph{Semidefinite programs and association
schemes}, Computing \textbf{63} (1999), 331--340.

\bibitem{indyk-2001} 
P. Indyk, \emph{Algorithmic applications of low-distortion geometric
embeddings}, In \emph{42nd Annual IEEE Symposium on Foundations of
Computer Science}, 10--33, 2001.

\bibitem{llr-1995}
N. Linial, E. London, and Y. Rabinovich, \emph{The geometry of graphs
and some of its algorithmic applications}, Combinatorica \textbf{15}
(1995), 215--246.

\bibitem{lm-2000}
N. Linial and A. Magen, \emph{Least-distortion Euclidean embeddings of
  graphs: Products of cycles and expanders}, J. Combin. Theory Ser. B
  \textbf{79} (1995), 157--171.

\bibitem{linial-2002}
N. Linial, \emph{Finite metric-spaces---combinatorics, geometry and
algorithms}, In \emph{Proceedings of the International Congress of
Mathematicians Vol. III (Beijing 2002)}, 573--586, 2002.

\bibitem{matousek-2002} 
J. Matousek, \emph{Lectures on discrete geometry}, Springer-Verlag,
New York, 2002.

\bibitem{vb-1996} 
L. Vandenberghe and S. Boyd, \emph{Semidefinite Programming}, SIAM
Rev. \textbf{38} (1996), 49--95.

\end{thebibliography}

\end{document}